\documentclass[graybox]{svmult}
\usepackage{fullwidth}
\usepackage{lipsum}  
\usepackage{url}
\usepackage{soul}
\usepackage{pslatex}
\usepackage{graphicx}
\usepackage{amsmath,amssymb,empheq,amsfonts,bm,marvosym,wasysym}
\usepackage{floatrow}
\usepackage{tikz}
\usepackage{pgfplots}
\usepackage{pgf}
\DeclareMathOperator{\EE}{\mathbb{E}}
\DeclareMathOperator{\V}{\mathbb{V}}
\newcommand\fracfun[1]{\textrm{frac}\left(#1\right)}
\newcommand\norm[1]{\left\lVert#1\right\rVert}

\newcommand\card[1]{|#1|}
\usepackage{graphicx}
\usepackage{amsmath}
\usepackage{amsfonts}
\usepackage{amssymb}
\usepackage{xcolor}
\usepackage{siunitx}
\usepackage{booktabs}
\usepackage{multirow}

\DeclarePairedDelimiter\floor{\lfloor}{\rfloor}
\usetikzlibrary{patterns,external,calc,angles,quotes}
\usepackage[]{algorithm2e}
\tikzexternaldisable
\usepackage{adjustbox}

\newcommand{\revAdd}[1]{{\color{black}{#1}}}

\usepackage{type1cm}        
\usepackage{makeidx}         
\usepackage{graphicx}        
\usepackage{multicol}        
\usepackage[bottom]{footmisc}

\usepackage{newtxtext}       %
\usepackage{newtxmath}       

\makeindex             

\begin{document}

\title*{On the Selection of Random Field Evaluation Points in the p-MLQMC Method}

\author{P. Blondeel, P. Robbe, S. Fran\c{c}ois, G. Lombaert and S. Vandewalle}
\institute{Philippe Blondeel, Pieterjan Robbe, Stefan Vandewalle, \at KU Leuven, Department of Computer Science
Celestijnenlaan 200A, 3001 Leuven, Belgium,\\ \email{{philippe.blondeel,pieterjan.robbe,stefan.vandewalle}@kuleuven.be}
\and Stijn Fran\c{c}ois, Geert Lombaert, \at KU Leuven, Department of Civil Engineering
Kasteelpark Arenberg 40, 3001 Leuven, Belgium,\\ \email{{stijn.francois,geert.lombaert}@kuleuven.be}
}
%
%
\maketitle

\abstract*{Engineering problems are often characterized by significant uncertainty in their material parameters. A typical example coming from geotechnical engineering is the slope stability problem where the soil's cohesion is modeled as a random field.  An efficient manner to account for this uncertainty is the novel sampling method called p-refined Multilevel Quasi-Monte Carlo (p-MLQMC). The p-MLQMC method uses a hierarchy of p-refined Finite Element meshes combined with a deterministic Quasi-Monte Carlo sampling rule. This combination yields a significant computational cost reduction with respect to classic Multilevel Monte Carlo. However, in previous work, not enough consideration was given \revAdd{to} how  to incorporate the uncertainty, modeled as a random field, in the Finite Element model with the p-MLQMC method. In the present work we investigate how this can be adequately achieved by means of the integration point method. We therefore investigate how the evaluation points of the random field are to be selected in order to obtain a variance reduction over the levels. We consider three different approaches. These approaches will be benchmarked on a slope stability problem in terms of computational runtime. We find that \revAdd{for a given tolerance} the \emph{Local Nested Approach} yields a speedup up to a factor five with respect to the \emph{Non-Nested approach}.}
\vspace{-0.5cm}
\abstract{Engineering problems are often characterized by significant uncertainty in their material parameters. A typical example coming from geotechnical engineering is the slope stability problem where the soil's cohesion is modeled as a random field.  An efficient manner to account for this uncertainty is the novel sampling method called p-refined Multilevel Quasi-Monte Carlo (p-MLQMC). The p-MLQMC method uses a hierarchy of p-refined Finite Element meshes combined with a deterministic Quasi-Monte Carlo sampling rule. This combination yields a significant computational cost reduction with respect to classic Multilevel Monte Carlo. However, in previous work, not enough consideration was given \revAdd{to} how  to incorporate the uncertainty, modeled as a random field, in the Finite Element model with the p-MLQMC method. In the present work we investigate how this can be adequately achieved by means of the integration point method. We therefore investigate how the evaluation points of the random field are to be selected in order to obtain a variance reduction over the levels. We consider three different approaches. These approaches will be benchmarked on a slope stability problem in terms of computational runtime. We find that \revAdd{for a given tolerance} the \emph{Local Nested Approach} yields a speedup up to a factor five with respect to the \emph{Non-Nested approach}.}
\vspace{-0.5cm}
\section{INTRODUCTION}

Starting from the work by Giles, see \cite{Giles,Giles2,Giles3},  we developed a novel multilevel method called p-refined Multilevel Quasi Monte Carlo (p-MLQMC), see \cite{Blondeel2020}. Similar to classic Multilevel Monte Carlo, see \cite{Giles}, p-MLQMC  uses a hierarchy of increasing resolution Finite Element meshes to achieve a computational speedup. Most of the samples are taken on coarse and computationally cheap meshes, while a decreasing number of samples are taken on finer and computationally expensive meshes. The major difference between classic Multilevel Monte Carlo and p-MLQMC resides in the refinement scheme used for constructing the mesh  hierarchy. In classic Multilevel Monte Carlo (h-MLMC), an h-refinement scheme is used to build the mesh hierarchy, see, for example \cite{Cliffe}. The accuracy of the model is increased by increasing the number of elements in the Finite Element mesh. In p-MLQMC, a p-refinement scheme is used to construct the mesh hierarchy. The accuracy of the model is increased by increasing the polynomial order of the element's shape functions while retaining the same number of elements. This approach reduces the computational cost with respect to h-MLMC, as shown in \cite{Blondeel2020}. Furthermore, instead of using a random sampling rule, i.e, the Monte Carlo method, p-MLQMC  uses a deterministic Quasi-Monte Carlo (QMC) sampling rule, yielding a further computational gain.

However, the p-MLQMC method presents the practitioner with a challenge. This challenge consists of adequately incorporating the uncertainty, modeled as a random field, into the Finite Element model. For classic Multilevel (Quasi)-Monte Carlo (h-ML(Q)MC) this is typically achieved by means of the midpoint method \cite{ChunChing}. The model uncertainty is represented as  scalars resulting from the evaluation of the random field at centroids of the elements. These scalars are then assigned to the elements. With this method, the uncertainty is modeled as being constant inside each element.  In h-refined multilevel methods, the midpoint method intrinsically links the spatial resolution of the mesh with the spatial resolution of the random field. An h-refinement of the mesh will result in a finer representation of the random field.
However, for the  p-MLQMC method, the midpoint method cannot be used. This is because the refinement scheme used in p-MLQMC does not increase the number of elements. The p-MLQMC method makes use of the integration point method, see \cite{MATTHIES1997283}. Scalars resulting from the evaluation of the random field at certain spatial locations are taken into account during numerical integration of the element stiffness matrices. With this method, the uncertainty varies inside each element. In the present work, we investigate how to adequately select  the spatial locations used for the evaluation of the random field. Specifically, we distinguish three different approaches \revAdd{to} how to select these random field evaluation points. The \emph{Non-Nested Approach} (NNA), the \emph{Global Nested Approach} (GNA), and the \emph{Local Nested Approach} (LNA). We investigate how these approaches affect the variance reduction in the p-MLQMC method, and how the total computational runtime increases over the levels. These approaches will be benchmarked on a model problem which consists of a slope stability problem, which assesses the stability of natural or man made slopes. The uncertainty is located in the soil's cohesion, and is represented as a two-dimensional lognormal random field.

The paper is structured as follows. First we give a theoretical background motivating our research, and give a concise overview of the building blocks of p-MLQMC. Second, we present the three approaches. Hereafter, we shortly discuss the underlying Finite Element solver, and introduce the model problem. Last, we present the results obtained with p-MLQMC for the three different approaches. Here we focus on the variance reduction over the levels and the effect on the total computational runtime.
\vspace{-0.5cm}
\section{Theoretical Background}
\label{sec:background}

Multilevel Monte Carlo methods rely on a hierarchy of meshes in order to achieve a speedup with respect to Monte Carlo. This speedup is achieved by \revAdd{writing the expected value of a quantity of interest on a fine mesh as the expected value of a quantity of interest on a coarse mesh together with a series of correction terms that express the difference in expected value of the quantity of interest on two successive finer meshes. In particular, given the hierarchy of approximations $P_0, P_1, \ldots, P_L$ for the quantity of interest $P$ computed on an increasingly finer mesh, we have the telescopic sum identity,

\begin{equation}
\EE[P_\text{L}]  =\EE[P_0]+\sum_{\ell=1}^\text{L} \EE[P_\ell -P_{\ell-1}].
\label{eq:telescoping}
\end{equation}}

\begin{figure}[H]
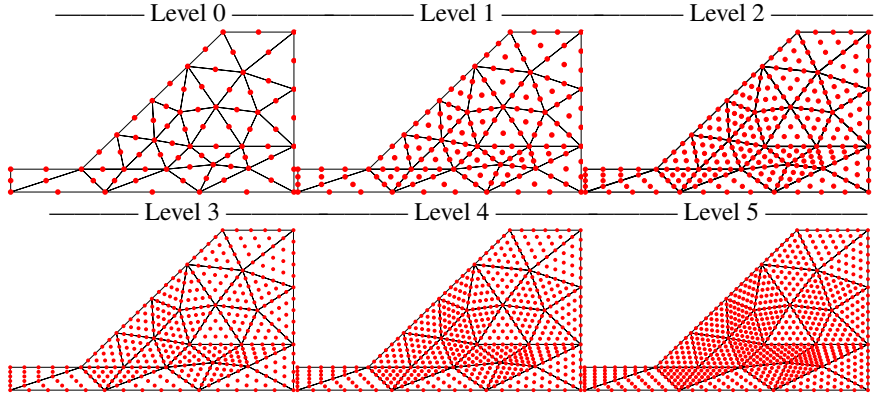

\hspace{0.08cm}
\begin{minipage}{0.3\textwidth}
\scalebox{0.8}{
\input{level0.tex}}
\end{minipage}
\begin{minipage}{0.3\textwidth}
\scalebox{0.8}{
\input{level1.tex}}
\end{minipage}
\begin{minipage}{0.3\textwidth}
\scalebox{0.8}{
\input{level2.tex}}
\end{minipage}
\\
\begin{minipage}{0.32\textwidth}
\scalebox{0.1}{
\input{L1.tex}}
\end{minipage}
\begin{minipage}{0.32\textwidth}
\scalebox{0.1}{
\input{L2.tex}}
\end{minipage}
\begin{minipage}{0.33\textwidth}
\scalebox{0.1}{
\input{L3.tex}}
\end{minipage}
\\
\begin{minipage}{0.3\textwidth}
\scalebox{0.8}{
\input{level3.tex}}
\end{minipage}
\begin{minipage}{0.3\textwidth}
\scalebox{0.8}{
\input{level4.tex}}
\end{minipage}
\begin{minipage}{0.3\textwidth}
\scalebox{0.8}{
\input{level5.tex}}
\end{minipage}
\\
\begin{minipage}{0.32\textwidth}
\scalebox{0.1}{
\input{L4.tex}}
\end{minipage}
\begin{minipage}{0.32\textwidth}
\scalebox{0.1}{
\input{L5.tex}}
\end{minipage}
\begin{minipage}{0.33\textwidth}
\scalebox{0.1}{
\input{L6.tex}}
\end{minipage}
\\
\vspace{-0.5cm}
\caption{p-refined hierarchy of \revAdd{approximations} used for the slope stability problem.}
\label{fig:meshes_2_p}
\end{figure}

This hierarchy of meshes can be obtained by applying an h-refinement scheme or a p-refinement scheme to a coarse mesh model. We opt for a hierarchy based on p-refinement. The hierarchy  applied to a discretized model of the slope stability problem is shown in Fig.\,\ref{fig:meshes_2_p}. Here, the Finite Element nodal points are represented as red dots. A more thorough discussion of the slope stability problem and the underlying Finite Element model is given in \S\ref{sec:FEM}.

\vspace{-0.5cm}
In the Multilevel Monte Carlo setting, the meshes in the hierarchy are commonly referred to as `levels'. The coarsest mesh is referred to as level 0. Subsequent finer meshes are assigned the next cardinal number, e.g., level 1, level 2, $\ldots$

The number of samples to be taken on levels  greater than 0 $\left(\ell>0\right)$, is proportional to the sample variance of the difference, $\V\left[\Delta P_\ell\right]$ with $\Delta P_\ell = P_\ell - P_{\ell-1}$ and $P$ a chosen quantity of interest (QoI).  It is only  for determining the number  of samples on level 0 that the sample variance $\V\left[P_\ell\right]$ is used. In order to obtain a decreasing number of samples per increasing level, i.e., $N_0 > N_1 > \cdots > N_\text{L}$, it is necessary to have a variance reduction over the levels, i.e., $\V\left[\Delta P_1\right] > \V\left[\Delta P_2\right] > \cdots > \V\left[\Delta P_\text{L}\right]$, and an increasing cost `of one solve' per increasing level, i.e., $\mathcal{C}_0 < \mathcal{C}_1 < \cdots < \mathcal{C}_\text{L}$. This variance reduction is only obtained when a strong positive correlation is achieved between the results of two successive levels. \revAdd{We have that}
\begin{equation}\label{eq:cov}
\begin{split}
\V\left[\Delta P_\ell\right] &= \V\left[P_\ell - P_{\ell-1}\right]\\
&= \V\left[P_\ell\right] + \V\left[P_{\ell-1}\right] -2\text{cov}\left(P_\ell,P_{\ell-1}\right),
\end{split}
\end{equation}
where $\text{cov}\left(P_\ell,P_{\ell-1}\right) = \rho_{\ell,\ell-1}\sqrt{\V\left[P_\ell\right]\V\left[P_{\ell-1}\right]}$ is the covariance between $P_\ell$ and $P_{\ell-1}$ with $\rho_{\ell,\ell-1}$ the correlation coefficient. The value of $\text{cov}\left(P_\ell,P_{\ell-1}\right)$ must be larger than 0 to have a large variance reduction, and hence an efficient multilevel method. 

In our p-MLQMC algorithm applied to a slope stability problem, see \cite{Blondeel2020}, the model uncertainty representing the soil's cohesion is located in  the elastic constitutive matrix $\mathbf{D}$. It is taken into account at the locations of the quadrature points when computing the integral in the  element stiffness matrices $\mathbf{K^e}$, i.e.,
\begin{equation}\label{eq:K}
\mathbf{K^e} = \int_{\Omega_e} \mathbf{B}^\text{T} \mathbf{D} \mathbf{B} d \Omega_e.
\end{equation} 
This is calculated in practice as
\begin{equation}\label{eq:intpoint}
\mathbf{K^e} = \sum_{i=1}^{\card{\mathbf{q}}}\mathbf{B}_i^\text{T} \mathbf{D}_i\mathbf{B}_i \text{w}_i,
\end{equation} 
where the matrix $\mathbf{B}_i = \mathbf{B}(\mathbf{q}^i)$ contains the derivatives of the shape functions,  evaluated at the  quadrature points $\mathbf{q}^i$,  the matrix $\mathbf{D}_i = $ \revAdd{$\mathbf{D}\left(\mathbf{x}^i,\omega\right)$}
 contains the model uncertainty \revAdd{computed at point $\mathbf{x}^i$}, and $\text{w}_i$ are the quadrature weights. The set of  quadrature points $\mathbf{q}$ is expressed in a local coordinate system of the triangular reference element. The uncertainty in the matrix  \revAdd{$\mathbf{D}_i$}  is represented by a scalar originating from the evaluation of the random field  at a carefully chosen spatial location.  This approach is commonly referred to as the \emph{integration point method} \cite{MATTHIES1997283}. Note that here the uncertainty is not constant in an element, i.e.,  $\mathbf{D}_1 \neq \mathbf{D}_2 \neq \cdots \neq \mathbf{D}_k$ in Eq.\,\eqref{eq:intpoint}. 
The scalar used in the matrix  \revAdd{$\mathbf{D}_i$}  originates from the evaluation of the random field at spatial location  \revAdd{$\mathbf{x}^i \in \mathbf{x}$}, in a global coordinate system of the mesh, by means of the Karhunen-Lo\`eve (KL) expansion with stochastic dimension $s$, i.e.,
\begin{equation}
Z(\mathbf{x},\omega)=\overline{Z}(\mathbf{x})+\sum_{n=1}^{s}  \sqrt{\theta_n} \xi_n(\omega) b_n(\mathbf{x})\,,
\label{eq:KLExpansion}
\end{equation}
where $\overline{Z}(\mathbf{x})$ is the mean of the field and  $\xi_n(\omega)$ denote i.i.d.\,standard normal random variables.
The symbols $\theta_n$  and $b_n(\mathbf{x})$  denote  the eigenvalues and eigenfunctions respectively, which are the solutions of the eigenvalue problem $\int_D C(\mathbf{x},\mathbf{y})b_n({\mathbf{y}})\mathrm{d}\mathbf{y} = \theta_n b_n({\mathbf{x}})$ with a given  covariance kernel $C(\mathbf{x},\mathbf{y})$. Note that in order to represent the uncertainty of the soil's cohesion in the considered slope stability problem, we do not use $Z(\mathbf{x},\omega)$ but $\text{exp}(Z(\mathbf{x},\omega))$, see \S\ref{sec:FEM}.

Our goal is to select evaluation points for Eq.\,\eqref{eq:KLExpansion}, grouped in sets $\left\lbrace\mathbf{x}_\ell\right\rbrace_{\ell=0}^\text{L}$, in order to ensure a good correlation between $P_\ell$ and $P_{\ell-1}$, i.e., such that the covariance, $\text{cov}\left(P_\ell,P_{\ell-1}\right)$, is as large as possible,  see Eq.\,\eqref{eq:cov}. We distinguish three different approaches for selecting the evaluation  points on the different levels, the \emph{Non-Nested Approach} (NNA), the \emph{Global Nested Approach} (GNA) and the \emph{Local Nested Approach} (LNA). All the approaches  start from the given sets of quadrature points on the different levels $\left\lbrace\mathbf{q}_\ell\right\rbrace_{\ell=0}^\text{L}$. Note that the number of quadrature points per level increases, $\card{\mathbf{q}_0}<\card{\mathbf{q}_1}\cdots<\card{\mathbf{q}_\text{L}}$. Given the sets $\left\lbrace\mathbf{q}_\ell\right\rbrace_{\ell=0}^\text{L}$, we select evaluation points for the random field in a local coordinate system and group them in sets $\left\lbrace\mathbf{x}^\text{local}_\ell\right\rbrace_{\ell=0}^\text{L}$, with the condition that $ \card{\mathbf{x}_\ell^{\text{local}}} = \card{\mathbf{q}_\ell}$. The points in the sets  $\left\lbrace\mathbf{x}^\text{local}_\ell\right\rbrace_{\ell=0}^\text{L}$ are then transformed to points in global coordinates, resulting in sets  $\left\lbrace\mathbf{x}_\ell\right\rbrace_{\ell=0}^\text{L}$. The points belonging to $\left\lbrace\mathbf{x}_\ell\right\rbrace_{\ell=0}^\text{L}$  are then used in Eq.\,\eqref{algo:NN}. Note that with the integration point method,  the spatial resolution of the field is  proportional to the  number of quadrature points. Increasing the number of quadrature points will result in a finer resolution of the random field.

Before elaborating further  upon these approaches, we first introduce the estimator used in our p-MLQMC algorithm. The estimator is given by
\begin{equation}
Q^{\textrm{MLQMC}}_\text{L}:= \frac{1}{R_0}\sum_{r=1}^{R_0}\frac{1}{N_0}\sum_{n=1}^{N_0} P_0(\mathbf{u}_0^{(r,n)}) + 
 \sum_{\ell=1}^\text{L} \frac{1}{R_\ell}\sum_{r=1}^{R_\ell}\left \{ \frac{1}{N_\ell} \sum_{n=1}^{N_\ell} \left( P_\ell(\mathbf{u}_\ell^{(r,n)})-P_{\ell-1}(\mathbf{u}_\ell^{(r,n)})\right) \right \}.
 \label{eq:MLQMC}
\end{equation}
It expresses the expected value of the quantity of interest on the finest level $\text{L}$  as the sample average of the quantity of interest on the coarsest level, plus a series of correction terms. 
\revAdd{A particularity of MLQMC consists of the use of deterministic sample points per level $\mathbf{u}_\ell^{\left(r,n\right)}$ in the unit cube, i.e., $[0,1]^s$, combined with an average  over a number of shifts $R_\ell$ on each level $\ell$. Averaging over the number of shifts is performed in order to obtain unbiased estimates of the computed stochastic quantities. 
The representation of these uniform distributed quasi-Monte Carlo sample points in $[0,1]^s$ is given by
\begin{equation}
\mathbf{u}^{\left(r,n\right)} = \fracfun{\phi_2(n)\mathbf{z} + \Xi_r} \quad \text{for} \quad n  \in \mathbb{N},
\label{eq:qmc_pt}
\end{equation}
where $\fracfun{x} = x - \floor{x}, x>0$,   $\phi_2$ is the radical inverse function in base 2, $\mathbf{z}$ is  an $s$-dimensional vector of positive integers,   $\Xi_r \in \left[0,1\right]^s$ is the random shift with $r=1,2,\ldots,R_\ell$, and $s$ the stochastic dimension. The representation of the points from Eq.\,\eqref{eq:qmc_pt} is known as a \emph{shifted rank-1 lattice rule}.  The generating vector $\mathbf{z}$ was constructed with the component-by-component (CBC) algorithm with decreasing weights, $\gamma_j = 1/j^2$, see \cite{KuoGenVec}. In the scope of this work, the uniform quasi-Monte Carlo sample points $\mathbf{u}_\ell^{\left(r,n\right)}$ are mapped from $[0,1]^s$ to $\mathbb{R}^s$ by means of the  inverse of the univariate standard normal cumulative distribution function,  ${\Phi}^{-1}(\cdot)$. The standard normal distributed quasi-Monte Carlo points, ${\Phi}^{-1}(\mathbf{u}^{\left(r,n\right)})$ are then substituted in Eq.\,\eqref{eq:KLExpansion}, and used as the random standard normal distributed $\xi_n(\omega)$ in order to generate  random field instances.}

\vspace{-0.7cm}
\section{Incorporating the uncertainty in the model}
\vspace{-0.2cm}
\label{sec:approaches}
In this section, we will discuss the mechanics behind the three approaches. We will show how the  evaluation points of the random field are selected on each level $\ell = \left\lbrace0,\ldots,\text{L}\right\rbrace$ for the different approaches.   Each of the approaches selects the evaluation points differently. However, all approaches start from the given set of the quadrature points $\mathbf{q}_\ell$. The points $\mathbf{q}^i_\ell \in \mathbf{q}_\ell$ are represented by $\textcolor{blue}{{\triangle}}_\ell$, on a reference triangular finite element on level $\ell$. Given the sets of quadrature points, the evaluation points of the random field, represented by $\textcolor{red}{\CIRCLE}_\ell$, are selected on a reference triangular finite element on level $\ell$, and grouped in the set $\mathbf{x}_\ell^\text{local}$.

The quadrature points  consist of a combination of points developed by Dunavant \cite{Dunavant} and Wandzurat  \cite{WANDZURAT20031829}, see Table\,\ref{Tab:References}. The code used to generate the Wandzurat points can be found at \cite{JohnBurkardt}.
\subsection{Non-Nested Approach}
\begin{figure}[H]
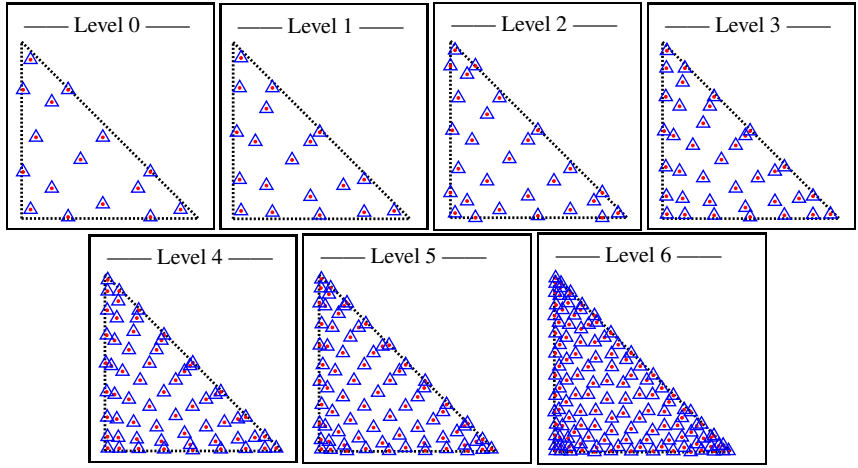

\begin{adjustbox}{varwidth=1.2\textwidth,fbox}
\begin{minipage}{0.18\textwidth}
\scalebox{0.7}{
\input{level00.tex}}
\end{minipage}

\begin{minipage}{0.18\textwidth}
\scalebox{0.26}{
\input{Triag_1_nn.tex}}
\end{minipage}
\end{adjustbox}
\begin{adjustbox}{varwidth=1.2\textwidth,fbox}
\begin{minipage}{0.18\textwidth}
\scalebox{0.7}{
\input{level11.tex}}
\end{minipage}

\begin{minipage}{.18\textwidth}
\scalebox{0.26}{
\input{Triag_2_nn.tex}}
\end{minipage}
\end{adjustbox}
\begin{adjustbox}{varwidth=1.2\textwidth,fbox}
\begin{minipage}{0.18\textwidth}
\scalebox{0.7}{
\input{level22.tex}}
\end{minipage}

\begin{minipage}{.18\textwidth}
\scalebox{0.26}{
\input{Triag_3_nn.tex}}
\end{minipage}
\end{adjustbox}
\begin{adjustbox}{varwidth=1.2\textwidth,fbox}
\begin{minipage}{0.18\textwidth}
\scalebox{0.7}{
\input{level33.tex}}
\end{minipage}

\begin{minipage}{.18\textwidth}
\scalebox{0.26}{
\input{Triag_4_nn.tex}}
\end{minipage}
\end{adjustbox}

\begin{adjustbox}{varwidth=1.2\textwidth,fbox}
\begin{minipage}{0.18\textwidth}
\scalebox{0.7}{
\input{level44.tex}}
\end{minipage}

\begin{minipage}{.18\textwidth}
\scalebox{0.26}{
\input{Triag_5_nn.tex}}
\end{minipage}
\end{adjustbox}
\vspace{0.5cm}
\begin{adjustbox}{varwidth=1.2\textwidth,fbox}
\begin{minipage}{0.2\textwidth}
\scalebox{0.7}{
\input{level55.tex}}
\end{minipage}

\begin{minipage}{0.2\textwidth}
\scalebox{0.26}{
\input{Triag_6_nn.tex}}
\end{minipage}
\end{adjustbox}
\begin{adjustbox}{varwidth=1.2\textwidth,fbox}
\begin{minipage}{0.2\textwidth}
\scalebox{0.7}{
\input{level66.tex}}
\end{minipage}

\begin{minipage}{.2\textwidth}
\scalebox{0.26}{
\input{Triag_7_nn.tex}}
\end{minipage}
\end{adjustbox}
\vspace{-0.8cm}
\caption{Locations of the quadrature points \textcolor{blue}{$\triangle$} and  of the evaluation points of the random field \textcolor{red}{$\CIRCLE$} on a reference triangular element in NNA.}
\label{fig:NonNested}
\end{figure}

For the Non-Nested Approach, the quadrature points on each level  are selected as the evaluation points of the random field, i.e., $\mathbf{x}^{\text{local}}_\ell = \mathbf{q}_\ell$ for  $\ell=\{0,\ldots,\text{L}\}$. Because the sets of quadrature points are not nested over the levels, i.e.,  $\mathbf{q}_0 \not\subseteq \mathbf{q}_1 \not\subseteq \cdots \not\subseteq \mathbf{q}_\text{L}$, it follows that the sets of the evaluation points of the random field are not nested, i.e., $\mathbf{x}^{\text{local}}_0 \not\subseteq \mathbf{x}^{\text{local}}_1 \not\subseteq \cdots \not\subseteq \mathbf{x}^{\text{local}}_\text{L}$, and thus $\mathbf{x}_0 \not\subseteq \mathbf{x}_1 \not\subseteq \cdots \not\subseteq \mathbf{x}_\text{L}$. This approach is the most straightforward one, and is illustrated in Fig.\,\ref{fig:NonNested}.  In Algorithm\,\ref{algo:NN}, we present the procedure which selects the  evaluation points of the random field for each level in local coordinates and groups them in sets $\lbrace \mathbf{x}^\text{local}_\ell\rbrace_{\ell=0}^\text{L}$.

\scalebox{0.9}{
\begin{algorithm}[H]
 \KwData{\\ Max level $\text{L}$,  Set of quadrature points per level $\left\lbrace \mathbf{q}_\ell \right\rbrace_{\ell=0}^\text{L}$}

$\ell \gets  L$\;
\While{$\ell \geq 0$}{
$\mathbf{x}^{\text{local}}_\ell \leftarrow \mathbf{q}_\ell$\;
$\ell \leftarrow \ell - 1$ \;
}
return $\left\lbrace \mathbf{x}^{\text{local}}_\ell \right\rbrace_{\ell=0}^\text{L}$
\vspace{0.4cm}
 \caption{Generation of the evaluation points of the random field in NNA.}
 \label{algo:NN}
\end{algorithm}
}

\subsection{Global Nested Approach}
\label{sec:gna}
\begin{figure}
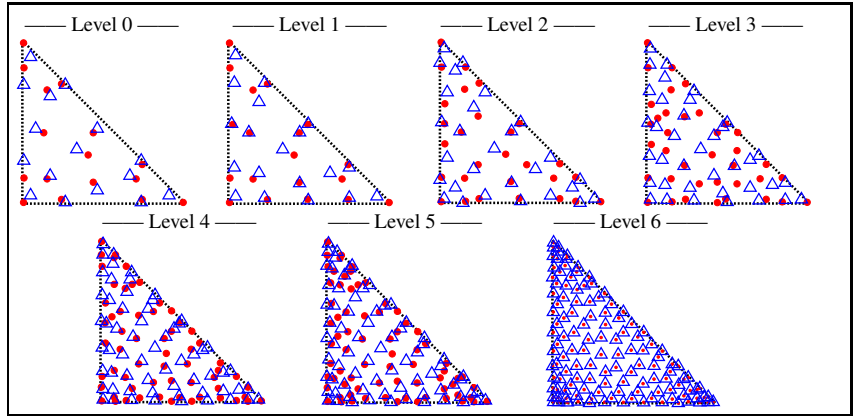

\begin{adjustbox}{varwidth=1.0\textwidth,fbox,center}
\begin{minipage}{0.23\textwidth}
\scalebox{0.65}{
\input{level00.tex}}
\end{minipage}
\begin{minipage}{0.23\textwidth}
\scalebox{0.65}{
\input{level11.tex}}
\end{minipage}
\begin{minipage}{0.23\textwidth}
\scalebox{0.65}{
\input{level22.tex}}
\end{minipage}
\begin{minipage}{0.23\textwidth}
\scalebox{0.65}{
\input{level33.tex}}
\end{minipage}
%

\begin{minipage}{0.23\textwidth}
\scalebox{0.24}{
\input{Triag_1_Glob.tex}}
\end{minipage}
\begin{minipage}{.23\textwidth}
\scalebox{0.24}{
\input{Triag_2_Glob.tex}}
\end{minipage}
\begin{minipage}{.23\textwidth}
\scalebox{0.24}{
\input{Triag_3_Glob.tex}}
\end{minipage}
\begin{minipage}{.23\textwidth}
\scalebox{0.24}{
\input{Triag_4_Glob.tex}}
\end{minipage}

\centering
\begin{minipage}{0.25\textwidth}
\scalebox{0.65}{
\input{level44.tex}}
\end{minipage}
\begin{minipage}{0.25\textwidth}
\scalebox{0.65}{
\input{level55.tex}}
\end{minipage}
\begin{minipage}{0.25\textwidth}
\scalebox{0.65}{
\input{level66.tex}}
\end{minipage}

\begin{minipage}{.25\textwidth}
\scalebox{0.24}{
\input{Triag_5_Glob.tex}}
\end{minipage}
\begin{minipage}{0.25\textwidth}
\scalebox{0.24}{
\input{Triag_6_Glob.tex}}
\end{minipage}
\begin{minipage}{.25\textwidth}
\scalebox{0.24}{
\input{Triag_7_Glob.tex}}
\end{minipage}
\end{adjustbox}
\vspace{-0.6cm}
\caption{Locations of the quadrature points \textcolor{blue}{$\triangle$} and  of the evaluation points of the random field \textcolor{red}{$\CIRCLE$} on a reference triangular element in GNA.}
\label{fig:GlobalNested}
\end{figure}

For the Global Nested Approach, we proceed in a different way. All the levels are correlated with each other. The sets of evaluation points of the random field are chosen  such that they are nested over all the levels, i.e., $\mathbf{x}^{\text{local}}_0 \subseteq \mathbf{x}^{\text{local}}_1 \subseteq \cdots \subseteq \mathbf{x}^{\text{local}}_\text{L}$, and thus $\mathbf{x}_0 \subseteq \mathbf{x}_1 \subseteq \cdots \subseteq \mathbf{x}_\text{L}$. For GNA, the sets of evaluation points are not equal to the sets of quadrature points, except on the finest level, i.e., $\mathbf{x}^{\text{local}}_\ell \neq \mathbf{q}_\ell$ $\text{for}$  $\ell=\{0,\ldots,\text{L}-1\}$ and $\mathbf{x}^{\text{local}}_\text{L} = \mathbf{q}_\text{L}$. The approach for selecting the evaluation points of the random field is as follows. The quadrature points on the finest level $\text{L}$ are selected as the evaluation points of the random field, i.e., $\mathbf{x}^{\text{local}}_\text{L} = \mathbf{q}_\text{L}$. The points selected for $\mathbf{x}^{\text{local}}_\ell$ on levels $\ell=\{\text{L}-1, \ldots, 0 \} $, consist of a number of points $\card{\mathbf{q}_\ell}$, which are selected from the set $\mathbf{x}^{\text{local}}_{\ell+1}$, such that each selected point is the closest neighbor of a point of the set $\mathbf{q}_\ell$, i.e., $\mathbf{x}_\ell^\text{local} :=  \underset{\underset{\card{\mathbf{x}_\ell^\text{local}}=\card{\mathbf{q}_\ell} }{\mathbf{x}_\ell^\text{local} \subseteq \mathbf{x}_{\ell+1}^\text{local}}}{\mathrm{argmin}} \textbf{D}\left(\mathbf{x}_{\ell+1}^\text{local},\mathbf{q}_\ell\right)$, where $\textbf{D}\left(\mathbf{a},\mathbf{b}\right):= \underset{a\in \mathbf{a}}{\sum}\textbf{d}\left(a,\mathbf{b}\right)$ is the distance between two sets, and where $\textbf{d}\left(a,\textbf{b}\right):= \text{inf}\left\lbrace\text{d}\left(a,b\right)|b\in\mathbf{b}\right\rbrace$ is the minimal distance between a point and a set, with $\text{d}\left(a,b\right)$ the Euclidean distance between two points. 
This is illustrated in Fig.\,\ref{fig:GlobalNested}. The procedure used to select  the evaluation points for GNA is given in Algorithm\,\ref{algo:Global}.

\scalebox{0.9}{
\begin{algorithm}[H]
 \KwData{\\ Max level $\text{L}$, Set of quadrature points per  level $\left\lbrace \mathbf{q}_\ell \right\rbrace_{\ell=0}^\text{L}$  }
$\mathbf{x}^{\text{local}}_\text{L} \leftarrow \mathbf{q}_\text{L}$\;
$\ell \gets  \text{L}-1$\;
\While{$\ell \geq 0$}{
$i \gets 1 $ \;
$\mathbf{x}^{\text{local}}_\ell \gets \varnothing$ \;
\While{$i \leq \card{\mathbf{q}_\ell} $}{ 
Find the point $p \in \mathbf{x}^{\text{local}}_{\ell+1}$ ,which is not in $\mathbf{x}^{\text{local}}_\ell$, closest to  $\mathbf{q}^i_\ell$   \;
$\mathbf{x}^{\text{local}}_\ell \leftarrow \mathbf{x}^{\text{local}}_\ell \cup \{p\}$; // Add it to the array \\
$i \leftarrow i + 1$ \;
}
$\ell \leftarrow \ell - 1$ \;
}

return $\left\lbrace \mathbf{x}^{\text{local}}_\ell \right\rbrace_{\ell=0}^\text{L}$
\vspace{0.4cm}
 \caption{Generation of the evaluation points of the random field in GNA.}
 \label{algo:Global}
\end{algorithm}}

\subsection{Local Nested Approach}
\begin{figure}[H]
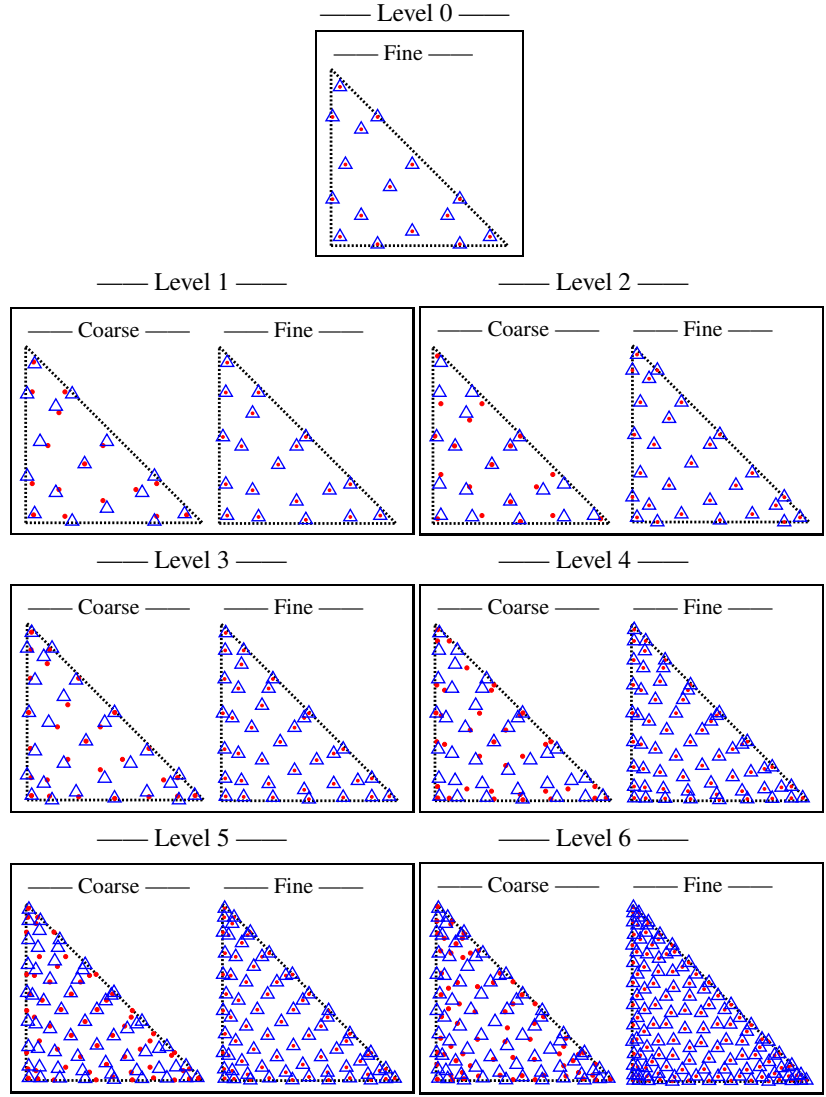

\begin{minipage}{0.25\textwidth}
\scalebox{0.8}{
\input{level00.tex}}
\end{minipage}

\begin{adjustbox}{varwidth=0.90\textwidth,fbox}
\begin{minipage}{0.24\textwidth}
\scalebox{0.7}{
\input{level_fine.tex}}
\end{minipage}

\begin{minipage}{0.24\textwidth}
\scalebox{0.26}{
\input{Triag_1_nn.tex}}
\end{minipage}
\end{adjustbox}

\vspace{0.1cm}
\begin{minipage}{0.45\textwidth}
\scalebox{0.8}{
\input{level11.tex}}
\end{minipage}
\vspace{0.1cm}
\begin{minipage}{0.3\textwidth}
\scalebox{0.8}{
\input{level22.tex}}
\end{minipage}

\begin{adjustbox}{varwidth=0.90\textwidth,fbox}
\begin{minipage}{0.24\textwidth}
\scalebox{0.7}{
\input{level_coarse.tex}}
\end{minipage}
\begin{minipage}{0.24\textwidth}
\scalebox{0.7}{
\input{level_fine.tex}}
\end{minipage}

\begin{minipage}{0.24\textwidth}
\scalebox{0.26}{
\input{Triag_21_Loc.tex}}
\end{minipage}
\begin{minipage}{.24\textwidth}
\scalebox{0.26}{
\input{Triag_22_Loc.tex}}
\end{minipage}
\end{adjustbox}
\begin{adjustbox}{varwidth=0.90\textwidth,fbox}
\begin{minipage}{0.24\textwidth}
\scalebox{0.7}{
\input{level_coarse.tex}}
\end{minipage}
\begin{minipage}{0.24\textwidth}
\scalebox{0.7}{
\input{level_fine.tex}}
\end{minipage}

\begin{minipage}{0.24\textwidth}
\scalebox{0.26}{
\input{Triag_32_Loc.tex}}
\end{minipage}
\begin{minipage}{.24\textwidth}
\scalebox{0.26}{
\input{Triag_33_Loc.tex}}
\end{minipage}
\end{adjustbox}
\vspace{0.1cm}

\begin{minipage}{0.45\textwidth}
\scalebox{0.8}{
\input{level33.tex}}
\end{minipage}
\vspace{0.1cm}
\begin{minipage}{0.3\textwidth}
\scalebox{0.8}{
\input{level44.tex}}
\end{minipage}
\begin{adjustbox}{varwidth=0.90\textwidth,fbox}
\begin{minipage}{0.24\textwidth}
\scalebox{0.7}{
\input{level_coarse.tex}}
\end{minipage}
\begin{minipage}{0.24\textwidth}
\scalebox{0.7}{
\input{level_fine.tex}}
\end{minipage}

\begin{minipage}{0.24\textwidth}
\scalebox{0.26}{
\input{Triag_43_Loc.tex}}
\end{minipage}
\begin{minipage}{.24\textwidth}
\scalebox{0.26}{
\input{Triag_44_Loc.tex}}
\end{minipage}
\end{adjustbox}
\begin{adjustbox}{varwidth=0.90\textwidth,fbox}
\begin{minipage}{0.24\textwidth}
\scalebox{0.7}{
\input{level_coarse.tex}}
\end{minipage}
\begin{minipage}{0.24\textwidth}
\scalebox{0.7}{
\input{level_fine.tex}}
\end{minipage}

\begin{minipage}{0.24\textwidth}
\scalebox{0.26}{
\input{Triag_54_Loc.tex}}
\end{minipage}
\begin{minipage}{.24\textwidth}
\scalebox{0.26}{
\input{Triag_55_Loc.tex}}
\end{minipage}
\end{adjustbox}
\vspace{0.1cm}

\begin{minipage}{0.45\textwidth}
\scalebox{0.8}{
\input{level55.tex}}
\end{minipage}
\vspace{0.1cm}
\begin{minipage}{0.3\textwidth}
\scalebox{0.8}{
\input{level66.tex}}
\end{minipage}

\begin{adjustbox}{varwidth=0.90\textwidth,fbox}
\begin{minipage}{0.24\textwidth}
\scalebox{0.7}{
\input{level_coarse.tex}}
\end{minipage}
\begin{minipage}{0.24\textwidth}
\scalebox{0.7}{
\input{level_fine.tex}}
\end{minipage}

\begin{minipage}{0.24\textwidth}
\scalebox{0.26}{
\input{Triag_65_Loc.tex}}
\end{minipage}
\begin{minipage}{.24\textwidth}
\scalebox{0.26}{
\input{Triag_66_Loc.tex}}
\end{minipage}
\end{adjustbox}
\begin{adjustbox}{varwidth=0.90\textwidth,fbox}
\begin{minipage}{0.24\textwidth}
\scalebox{0.7}{
\input{level_coarse.tex}}
\end{minipage}
\begin{minipage}{0.24\textwidth}
\scalebox{0.7}{
\input{level_fine.tex}}
\end{minipage}

\begin{minipage}{0.24\textwidth}
\scalebox{0.26}{
\input{Triag_76_Loc.tex}}
\end{minipage}
\begin{minipage}{.24\textwidth}
\scalebox{0.26}{
\input{Triag_77_Loc.tex}}
\end{minipage}
\end{adjustbox}
\vspace{-0.5cm}
\caption{Locations of the quadrature points \textcolor{blue}{$\triangle$} and  of the evaluation points of the random field \textcolor{red}{$\CIRCLE$} on a reference triangular element in LNA.}
\label{fig:Local}
\end{figure}

\revAdd{As was the case for the previous approaches, the user first defines which sets of quadrature points $\mathbf{q}_\ell$ are to be used. Here we set $\mathbf{q}_{\ell,\text{fine}} := \mathbf{q}_\ell$, and $\mathbf{q}_{\ell,\text{coarse}} := \mathbf{q}_{\ell-1,\text{fine}}$.}
The Local Nested Approach is as follows.  Rather than correlating all the levels with each other, we now correlate them two-by-two. Each level $\ell=\{1,\ldots,\text{L}\}$ has two sets of evaluation points $\mathbf{x}_{\ell,\text{coarse}}$ and $\mathbf{x}_{\ell,\text{fine}}$, which are nested, i.e., $\mathbf{x}_{\ell,\text{coarse}} \subseteq \mathbf{x}_{\ell,\text{fine}}$. The  points in these sets are used to generate a coarse and a fine representation of the random field on level $\ell$. NNA and GNA  have only one set of points per level, and thus only one representation of the random field per level. The coarse representation of the random field essentially acts as a representation of the field on level $\ell-1$. This is because  $\mathbf{q}_{\ell,\text{coarse}} = \mathbf{q}_{\ell-1,\text{fine}}$. The selection process is as follows. For each  level $\ell =\left\lbrace0,\ldots,\text{L}\right\rbrace$, $\mathbf{x}^{\text{local}}_{\ell,\text{fine}} = \mathbf{q}_{\ell,\text{fine}}$. The points in $\mathbf{x}^{\text{local}}_{\ell,\text{coarse}}$ are selected according to the same methodology as in GNA, i.e., they are selected from the set $\mathbf{x}^{\text{local}}_{\ell,\text{fine}}$, such that each selected point is the closest neighbor to a point of the set $\mathbf{q}_{\ell,\text{coarse}}$. This is illustrated in  Fig.\,\ref{fig:Local}. The main advantage of this approach is level exchangeability and extensibility. With exchangeability we mean that if one pair of correlated levels, say $\tau$ and $\tau-1$, exhibits a `sub-optimal' value of $\V\left[\Delta P_\tau\right]$ with respect to the variances $\V\left[\Delta P_\ell\right]$ on other levels, this pair can easily be exchanged against another newly computed pair with a different set of quadrature points. This is in contrast with GNA where the whole hierarchy needs to be recomputed. With level extensibility we mean that if for a user requested tolerance $\varepsilon$ and  maximum level $\text{L}$, the tolerance is not reached, the hierarchy can easily be extended  by supplying the extra needed level(s) and reusing the previously computed samples. In case of GNA, the whole hierarchy needs to be recomputed with extra level(s) and the previously computed samples cannot be reused. This level extensibility is the major advantage of LNA over GNA.

\revAdd{An important note must be made concerning the LNA approach. While it successfully correlates the solutions of two successive levels, the expected value obtained from the telescoping sum is biased. We have observed a small bias of the order of $10^{-6}$ with respect to the actual values, an error that is well below the discretization error of the finite element discretization. The reasons behind this additional bias stems from the fact that substitute random fields are used. We are currently investigating how this additional bias can be avoided.}

\revAdd{
The procedure used to generate the point set for LNA is given in Algorithm\,\ref{algo:Local}. For LNA, each level has two representations of the random field, a coarse and a fine one, with sets  $\mathbf{x}_{\ell,\text{coarse}}$ and $\mathbf{x}_{\ell,\text{fine}}$.
\begin{fullwidth}[width=\linewidth+3cm,leftmargin=-1.5cm,rightmargin=-1.5cm]
\vspace{2pt}
\scalebox{0.8}{
\begin{algorithm}[H]
 \KwData{\\ Max level $\text{L}$, Set of quadrature points per  level $\left\lbrace \mathbf{q}_{\ell,\text{coarse}} \right\rbrace_{\ell=0}^\text{L}$ and $\left\lbrace \mathbf{q}_{\ell,\text{fine}} \right\rbrace_{\ell=0}^\text{L}$}
$\ell \gets \text{L}$ \;

\While{$\mathcal{\ell}>0$}{
$\mathbf{x}^{\text{local}}_{\ell,\text{fine}} \gets \mathbf{q}_{\ell,\text{fine}}$\;
$\mathbf{x}^{\text{local}}_{\ell,\text{coarse}} \gets \varnothing$ \;
$i \gets 1 $ \;
$k \gets \card{\mathbf{q}_{\ell,\text{coarse}}}$\;
\While{$i \leq  k$}{ 
Find the point $p \in \mathbf{x}^{\text{local}}_{\ell,\text{fine}}$ ,which is not in $\mathbf{x}^{\text{local}}_{\ell,\text{coarse}}$, closest to  $\mathbf{q}^i_{\ell,\text{coarse}}$   \;
$\mathbf{x}^{\text{local}}_{\ell,\text{coarse}} \leftarrow \mathbf{x}^{\text{local}}_{\ell,\text{coarse}} \cup \{p\}$; // Add it to the array \\

$i \leftarrow i + 1$ \;
}
$\ell \gets  \ell-1$\;

}
$\mathbf{x}^{\text{local}}_{0,\text{fine}} \leftarrow \mathbf{q}_{0,\text{fine}}$\;
return $\left\lbrace \mathbf{x}^{\text{local}}_{0,\text{fine}}, \mathbf{x}^{\text{local}}_{1,\text{fine}}, \mathbf{x}^{\text{local}}_{1,\text{coarse}}, \mathbf{x}^{\text{local}}_{2,\text{fine}}, \mathbf{x}^{\text{local}}_{2,\text{coarse}}, \ldots \right\rbrace$
\vspace{0.4cm}
 \caption{Generation of the evaluation points of the random field in  LNA.}
 \label{algo:Local}
\end{algorithm}}
\vspace{-0.9cm}
\end{fullwidth}
}

\vspace{-1.3cm}

\revAdd{\subsection{Discussion of the Computational Cost}
\vspace{-0.4cm}
We will now discuss how the computational cost for each of the approaches can be determined. This is done regardless of the number of samples. The total computational cost is split in an offline part, i.e., the cost of computing the eigenvalues and eigenvectors of the random fields, and an online part, i.e., the cost of computing point evaluations of the random fields at the evaluation points. The offline cost is only accounted for at startup, while the online cost is accounted for when computing each sample. 
For NNA, the offline cost is equal to $\sum_{\ell=0}^\text{L} C^{\text{eig}}_\ell$, with $C^{\text{eig}}_\ell$ the cost of computing the eigenvalues and eigenvectors of the random field on level $\ell$. For LNA, this cost is equal to $\sum_{\ell=0}^\text{L} C^{\text{eig}}_{\ell,\text{fine}}$. Only the eigenvalues and eigenvectors of the fine representation of the fields on each level need be computed. Note that $C^{\text{eig}}_\ell = C^{\text{eig}}_{\ell,\text{fine}}$. For GNA, the offline cost equals $\sum_{\ell=0}^\text{L} C^{\text{eig}}_\ell$. In case of GNA, the choice could be made to compute  the eigenvectors and  eigenfunctions  only on the finest level L, resulting in the offline cost $C^{\text{eig}}_\text{L}$. This is only possible because of the property $\mathbf{x}_0 \subseteq \mathbf{x}_\ell \subseteq \cdots \subseteq \mathbf{x}_\text{L}$, see \S\,\ref{sec:gna}. In practice this will not be done, because of a drastic increase of the online cost, see further on. The online cost for NNA on levels $\ell > 0$ is equal to the cost of computing point evaluations of the random field on level $\ell$, $C^{\text{samp}}_\ell$, and on level $\ell-1$, $C^{\text{samp}}_{\ell-1}$. This yields a total cost per sample equal to  $C^{\text{samp}}_\ell+C^{\text{samp}}_{\ell-1}$. For LNA, the online cost on level $\ell > 0$ is only equal to the cost of computing point evaluations of the fine representation of the field on level $\ell$, $C^{\text{samp}}_{\ell,\text{fine}}$. In order to represent the coarse field on level $\ell$, a restriction of the point evaluations of the fine random field on level $\ell$ is taken. This is because of the property $\mathbf{x}_{\ell,\text{coarse}} \subseteq \mathbf{x}_{\ell,\text{fine}}$, see Fig.\,\ref{fig:Local}. There is no cost associated with the restriction. Note that $C_\ell^\text{samp}=C_{\ell,\text{fine}}^\text{samp}$. For GNA, the online cost on level $\ell>0$ amounts to $C^{\text{samp}}_\ell$, if the eigenfunctions and eigenvectors have been computed for each level $\ell$. Then, the online cost is equal to the one of LNA. However, if the eigenfunctions and eigenvalues have only been computed on the finest level, the online cost for GNA equals $C^{\text{samp}}_\text{L}$ regardless of the level. Point evaluations of the random field are computed on level L and restricted to the desired level $\ell$.  In practice this is not done because of the much higher cost this incurs, $C^{\text{samp}}_\text{L} \gg C^{\text{samp}}_\ell$.}
\vspace{-0.7cm}
\section{Model Problem}
\label{sec:FEM}
The model problem we consider for benchmarking the three approaches consists of a slope stability problem  where the soil's cohesion has a spatially varying uncertainty \cite{Whenham}. In a slope stability problem the safety of the slope can be assessed by evaluating the vertical displacement of the top of the slope when sustaining its own weight.  We consider the displacement in the plastic domain, which is governed by the Drucker--Prager yield criterion. A small amount of isotropic linear hardening is taken into account for numerical stability reasons. Because of the nonlinear stress-strain relation arising in the plastic domain, a Newton--Raphson iterative solver is used.  In order to compute the displacement in a slope stability problem, an incremental load approach is used, i.e., the total load resulting from the slope's weight is added in steps starting with a force of $0\,\mathrm{N}$ until the downward force resulting from the slope's weight is reached.  This approach results in the following  system of equations for the displacement,
\begin{equation}\label{Displacement_eq_plast}
\mathbf{K} \Delta\mathbf{u} = \mathbf{r},
\end{equation}
where  $\Delta\mathbf{u}$ stands for the displacement increment, $\mathbf{K}$ the global stiffness matrix resulting from the assembly of element stiffness matrices $\mathbf{K^e}$, see Eq.\,\eqref{eq:K}. The vector $\mathbf{r}$ is the residual, 
\begin{equation}
\mathbf{r}=\mathbf{f}+\Delta\mathbf{f}-\mathbf{k},
\end{equation}
where  $\mathbf{f}$ stands for the sum of the external force increments applied in the previous steps, $\Delta\mathbf{f}$  for the applied load increment of the current step and   $\mathbf{k}$  for the internal force resulting from the stresses. For a more thorough explanation on the methods used to solve the slope stability problem we refer to \cite[Chapter 2 $\S$4 and Chapter 7 $\S$3 and $\S$4]{Borst}. 

The mesh hierarchy shown in Fig.\,\ref{fig:meshes_2_p} is generated by using a  combination of the open source mesh generator GMSH \cite{GMSH} and \textsc{Matlab} \cite{MATLAB:2017}. Table\,\ref{Tab:References} lists the number of elements (Nel), degrees of freedom (DOF), element order  per level (Order),  the number of quadrature points per element (Nquad), and the reference for the quadrature points (Ref) per level for p-MLQMC. The number of quadrature points  is chosen as to increase the spatial resolution of the field per increasing level, and to ensure numerical stability of the computations of the displacement in the plastic domain. In this paper we consider two-dimensional uniform, Lagrange triangular elements. 
\begin{table}[H]
\caption{Number of elements, degrees of freedom, element order and number of quadrature points for the model problem.}
\label{Tab:References}
\centering
\scalebox{1.0}{
\begin{tabular}{cccccc}
\toprule
 \multicolumn{6}{c}{p-MLQMC} \\
\cmidrule(rl{4pt}){1-6}  
 {Level} &  {Nel} &   {DOF} & {Order} & {Nquad} & Ref.\\
 0 &  33 &  160 &2      &16 &\cite{Dunavant}  \\
 1 &  33 &  338 &3      &19 &\cite{Dunavant} \\
 2 &  33 &  582 &4      &28 &\cite{Dunavant}\\
 3 &  33 &   892 &5      &37 &\cite{Dunavant}\\ 
 4 &  33 &   1268  &6      &61 &\cite{Dunavant}\\ 
 5 &  33 &   1710  &7      &73 &\cite{Dunavant}\\ 
 6 &  33 &   2218  &8      &126 &\cite{WANDZURAT20031829}\\ 
\bottomrule
\end{tabular}}
\end{table} 
 
\begin{figure}[H]
\sidecaption
\scalebox{0.1}{
\input{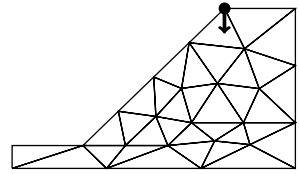}}
\caption{The vertical displacement of the upper left node as the QoI, indicated with an arrow.}
\label{fig:QoI2}
\end{figure}

We consider  the vertical displacement in meters of the upper left node of the model as a quantity of interest (QoI). This is depicted in Fig.\,\ref{fig:QoI2} by the arrow.
The uncertainty of the soil's cohesion is represented by means of  a  lognormal random field. This field is obtained by applying the exponential to the field obtained in Eq.\,\eqref{eq:KLExpansion}, $Z_{\text{lognormal}}(\mathbf{x},\omega)=\exp(Z(\mathbf{x},\omega))$. 
For the covariance \revAdd{k}ernel of the random field, we use the Mat\'ern covariance kernel,
\vspace{-0.1cm}
\begin{equation}\label{eq:covariance_kernel}
C(\mathbf{x},\mathbf{y}):= \sigma^2 \frac{1}{2^{\nu-1} \Gamma\left(\nu\right)}\left(\sqrt{2\nu}\dfrac{\norm{\mathbf{x}-\mathbf{y}}_2}{\lambda}\right)^\nu K_\nu \left(\sqrt{2\nu}\dfrac{\norm{\mathbf{x}-\mathbf{y}}_2}{\lambda}\right) \,,
\end{equation}
with $\nu = 2.0$ the smoothness parameter, $K_\nu$ the modified Bessel function of the second kind, $\sigma^2 = 1$ the variance and $\lambda = 0.3$ the correlation length. The characteristics of the lognormal distribution used to represent the uncertainty of the soil's cohesion are as follows: a mean of $8.02\,\mathrm{kPa}$ and a standard deviation of $400\,\mathrm{Pa}$. The spatial dimensions of the slope are: a length of $20\,\mathrm{m}$, a height of $14\,\mathrm{m}$ and a slope angle of $\ang{30}$. The material characteristics are: a Young's modulus of $30\,\mathrm{MPa}$, a Poisson ratio of $0.25$, a density of $1330\,\mathrm{kg/m^3}$ and a friction angle of $\ang{20}$. Plane strain is considered for this problem. The number of stochastic dimensions considered for the generation of the Gaussian random field is $s \!= 400$, see Eq.\,\eqref{eq:KLExpansion}. With a value $s \!= 400$, 99$\%$ of the variability of the random field is accounted for.

The stochastic part of our simulations was performed with the Julia packages \textbf{MultilevelEstimators.jl}, see \cite{PieterJanGit1}, and \textbf{GaussianRandomFields.jl}, see \cite{PieterJanGit2}. The Finite Element code used is an in-house \textsc{Matlab} code developed by the Structural Mechanics Section of the KU Leuven. All the results have been computed  on a  workstation equipped with 2 physical cores, Intel Xeon E5-2680 v3 CPU's, each with 12 logical cores, clocked at 2.50 GHz,  and a total of 128 GB RAM. 
\vspace{-0.6cm}
\section{Numerical Results}
\vspace{-0.4cm}
In this section we present our numerical results obtained with the p-MLQMC method.
\vspace{-1.2cm}
\subsection{Displacement}
\begin{figure}[H]
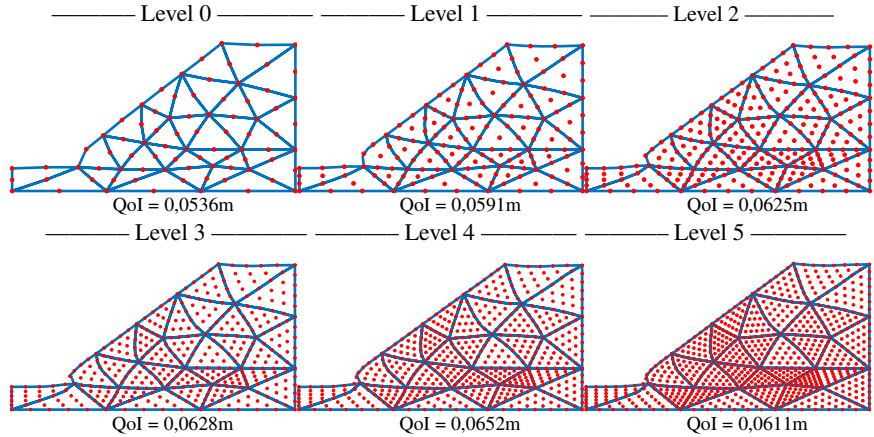

\hspace{0.08cm}
\begin{minipage}{0.3\textwidth}
\scalebox{0.75}{
\input{level0.tex}}
\end{minipage}
\begin{minipage}{0.3\textwidth}
\scalebox{0.75}{
\input{level1.tex}}
\end{minipage}
\begin{minipage}{0.3\textwidth}
\scalebox{0.7}{
\input{level2.tex}}
\end{minipage}
\\
\begin{minipage}{0.32\textwidth}
\scalebox{0.1}{
\input{Displacement/M0_S0}}
\end{minipage}
\begin{minipage}{0.32\textwidth}
\scalebox{0.1}{
\input{Displacement/M1_S1}}
\end{minipage}
\begin{minipage}{0.32\textwidth}
\scalebox{0.1}{
\input{Displacement/M2_S2}}
\end{minipage}
\\
\begin{minipage}{0.3\textwidth}
\scalebox{0.75}{
\input{level3.tex}}
\end{minipage}
\begin{minipage}{0.3\textwidth}
\scalebox{0.75}{
\input{level4.tex}}
\end{minipage}
\begin{minipage}{0.3\textwidth}
\scalebox{0.75}{
\input{level5.tex}}
\end{minipage}
\\
\begin{minipage}{0.32\textwidth}
\scalebox{0.1}{
\input{Displacement/M3_S3}}
\end{minipage}
\begin{minipage}{0.32\textwidth}
\scalebox{0.1}{
\input{Displacement/M4_S4}}
\end{minipage}
\begin{minipage}{0.32\textwidth}
\scalebox{0.1}{
\input{Displacement/M5_S5}}
\end{minipage}

\caption{ \revAdd{Displaced} mesh\revAdd{es} and QoI for different samples of the random field.}
\label{fig:Disp}
\end{figure}

\vspace{-0.2cm}
In Fig.\,\ref{fig:Disp}, we show the \revAdd{displaced}  mesh\revAdd{es} and the\revAdd{ir} nodes for a single sample of the random field on different levels. For better visualization, the displacement of the mesh and nodes in the figure have been exaggerated by a factor 20.  The value of the QoI is listed beneath each figure depicting the displacement.  
\vspace{-0.2cm}

\subsection{Variance and Expected value over the Levels}
\label{sec:var}
\begin{figure}[H]
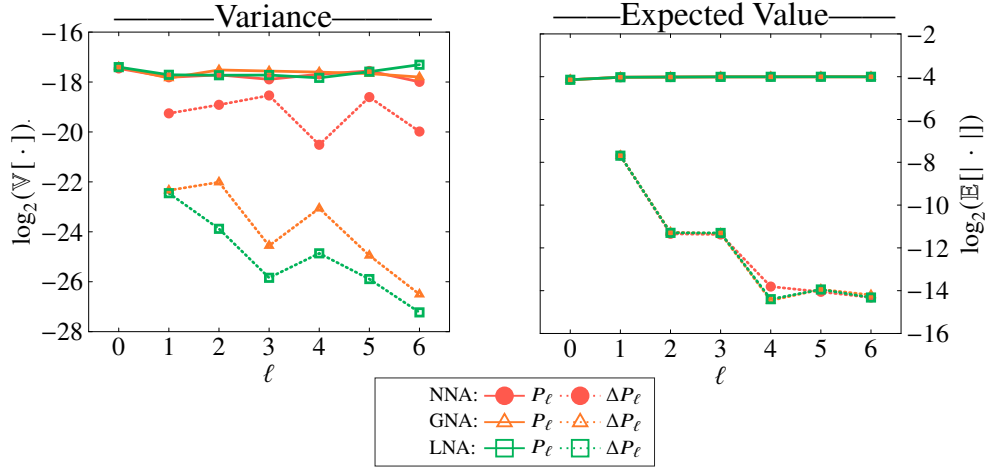
 
\hspace{-0.7cm}  
         \begin{minipage}{.6\textwidth}
        \scalebox{0.7}{
\input{Rates/V_case}}
 \end{minipage} 
 \begin{minipage}{.4\textwidth}
        \scalebox{0.7}{
\input{Rates/E_case}}
 \end{minipage}   
   
 \vspace{-0.3cm}
 \begin{minipage}{0.4\textwidth}
 \hspace{-1.3cm}
        \scalebox{0.9}{
\input{Rates/General_Legend}}
 \end{minipage}
 \vspace{-3.9cm}
  \caption{Variance and Expected Value over the levels.}
   \label{fig:variance}
\end{figure}

In Fig.\,\ref{fig:variance} we show the sample variance over the levels $\V\left[P_\ell\right]$, the sample variance of the difference over the levels $ \V\left[\Delta P_\ell\right]$, the expected value over the levels $\EE\left[P_\ell\right]$ and the expected value of the difference over the levels $ \EE\left[\Delta P_\ell\right]$.

\vspace{-0.5cm}
\vspace{-0.6cm}
As expected with multilevel methods, we observe that $\EE\left[P_\ell\right]$ remains constant over the levels, while $ \EE\left[\Delta P_\ell\right]$  decreases with increasing level. This is the case for all approaches we introduced in \S\ref{sec:approaches}.

As explained in \S\ref{sec:background}, multilevel methods are based on a variance reduction by means of a hierarchical refinement of Finite Element meshes. In practice this means that the sample variance  $\V\left[P_\ell\right]$ remains constant across the levels, while the sample variance of the difference over the levels $ \V\left[\Delta P_\ell\right]$ decreases for increasing level.
This is indeed what we observe for GNA and LNA.  For NNA we observe that $ \V\left[\Delta P_\ell\right]$ does not decrease. 
From Fig.\,\ref{fig:variance}, we can conclude that the choice of using nested over non-nested spatial locations over the levels as evaluation points for the random field  greatly improves the behavior of $ \V\left[\Delta P_\ell\right]$.  This influence stems from  a `bad' correlation between the results of two successive levels in the NNA case, see Eq.\,\eqref{eq:cov}. We will show in the next section that the number of samples per level required by NNA will be larger than GNA or LNA. This will impact the total runtime.
\subsection{Number of Samples}

\begin{figure}[H]
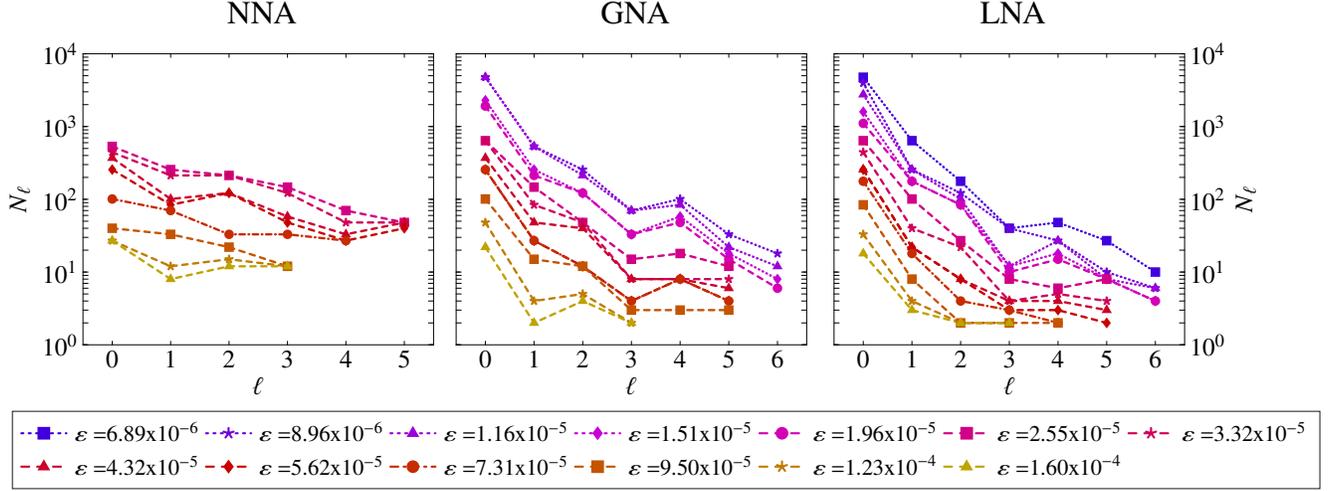

 \hspace{-3.cm}
 \begin{minipage}{0.3\textwidth}
                 \scalebox{.68}{
\input{samples/samples_1}
}
 \end{minipage}  
  \hspace{2.2cm}
  \begin{minipage}{0.3\textwidth}
                 \scalebox{.68}{
\input{samples/samples_2}
}
 \end{minipage} 
   \hspace{1.5cm}
  \begin{minipage}{0.3\textwidth}
                 \scalebox{.68}{
\input{samples/samples_3}
}
 \end{minipage} 
 
 \begin{minipage}[b]{1.\textwidth}
 \hspace{-2.8cm}
\scalebox{1.}{
\input{samples/General_Legend_case2}}
\end{minipage}
\vspace{-5.2cm}
 \caption{Number of samples for the three implementations for a given tolerance $\varepsilon$.}
 \label{fig:samples}
\end{figure}

In Fig.\,\ref{fig:samples}, we show the number of samples for the three approaches for thirteen different tolerances on the RMSE. These numbers do not include the number of shifts, which value is taken to be $R_\ell=10$ $\text{for}$ $\ell=\left\lbrace 0,\ldots,\text{L}\right\rbrace$.
We observe that for a given tolerance $\varepsilon$, the number of samples for NNA for levels greater than 0 is higher than for  GNA and  LNA. This is due to the slow decrease of $ \V\left[\Delta P_\ell\right]$, see Fig.\,\ref{fig:variance}.

\revAdd{
Unlike Multilevel Monte Carlo, the number of sample per level for MLQMC is not the result of an optimization problem,  see 
\cite{Giles}. The number of samples in MLQMC are chosen according to an adaptive `doubling' algorithm, where the number of samples is each time multiplied by a factor until the statistical constrain are satisfied. For a more thorough explanation we refer to \cite[\S 5]{Giles3}. The advantage of an adaptive algorithm consists of the fact that no equation expressing the  number of samples in function of the variance has to be derived. Indeed, such an expression would require the evaluation of terms which are difficult to estimate at runtime. In our approach, the sample multiplication factor is chosen as 1.2. The multiplication factor of 1.2 is chosen over the more naive value of 2 because, in the current setting  each sample involves the solution of a complex PDE, where the doubling of the number of samples would lead to huge jumps in the total cost of the estimator. The reduction from 2 to  1.2 is a compromise leading to a more gradual increase of the computational cost, while ensuring the enough progress is made.}

\subsection{Runtimes}
\begin{figure}[H]
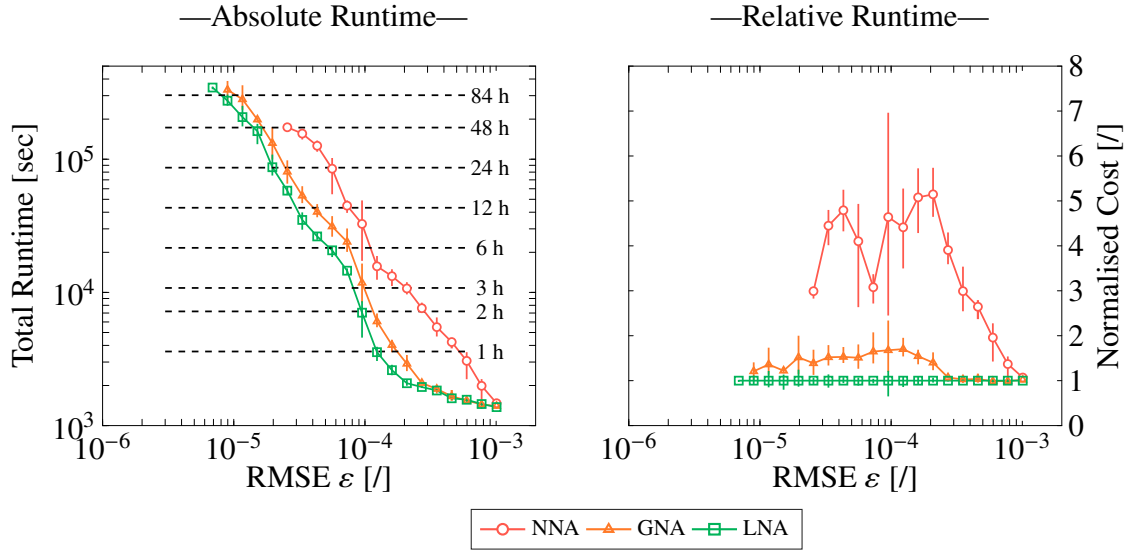

\hspace{-2.8cm}
 \begin{minipage}{0.45\textwidth}
                 \scalebox{.84}{
\input{Runtime/run_time_case_MEAN}
}
 \end{minipage}
 \hspace{2.4cm}
 \begin{minipage}{0.45\textwidth}
                 \scalebox{.84}{
\input{Runtime/run_time_case_Normalized_MEAN}
}
 \end{minipage}
     
 \begin{minipage}[b]{0.6\textwidth}
\scalebox{1.}{
\hspace{-0.5cm}
\input{Runtime/General_Legend_case2}}
\end{minipage}
 \vspace{-5.3cm}
\caption{Runtimes  \revAdd{as a} function of requested user tolerance.}
\label{fig:runtime_case}
\end{figure}

We show the absolute and relative runtime as a function of the user requested tolerance $\varepsilon$ on the RMSE for the different implementations in Fig.\,\ref{fig:runtime_case}.

\vspace{-.3cm}
\revAdd{For the shown results, the maximal number of available levels in the p-MLQMC method has been set to 7, i.e., $\text{L}$ = 6. The algorithm adaptively choses the number of levels needed to satisfy a given tolerance on the RMSE during runtime. If for a given tolerance on the RMSE, 8 levels are needed, an extra level can easily be added in case of the LNA approach while reusing the previously computed samples. For GNA, an extra level can also be added but the previously computed samples can not be reused in the new hierarchy of meshes. }

The results for the absolute runtime are expressed in seconds. For the relative runtime, we have normalized the computational cost of all three approaches such that the results for LNA for each tolerance have unity cost.  For both the results of the absolute and the relative runtime, we show the average, computed over three independent simulation runs, together with the minimum and maximum bounds.

We observe that \revAdd{for a given tolerance}, LNA achieves a speedup up to a factor 5 with respect to NNA and a factor 1.5 with respect to GNA. GNA also outperforms NNA in terms of lower computational cost.  This better performance of GNA and LNA is due to a lower number of samples per level, resulting from a better correlation between the successive levels. We can thus state that the Local and Global Nested Approaches achieve a lower computational cost with respect to the Non-Nested Approach \revAdd{for a given tolerance}.

\section{CONCLUSIONS}
In this work, we investigated how the spatial locations used for the evaluation of the random field by means of a Karhunen-Lo\`eve expansion impact the performance of the p-MLQMC method. We distinguished three different approaches, the \emph{Non-Nested Approach}, the \emph{Global Nested Approach} and the \emph{Local Nested Approach}.  We demonstrated that the choice of the evaluation points of the random field impacts the variance reduction over the levels $\V\left[\Delta P_\ell \right]$. We showed that the Global and Local Nested approaches exhibit a much better decrease of $\V\left[\Delta P_\ell \right]$ due to a better correlation between the levels than the Non-Nested Approach. This leads to  a lower number of samples for the  Nested Approaches and thus a lower total runtime for a given tolerance. Furthermore we have shown that the Local Nested Approach has the additional properties of level exchangeability and extensibility with respect to the Global Nested Approach. By correlating the levels two-by-two, in the Local Nested Approach, one pair of levels can easily be exchanged for another computed pair, if needed. In addition, the hierarchy can also easily be extended by adding a newly computed pair. When exchanging a pair of levels or extending the hierarchy, the previously computed samples can be reused. The Global Nested Approach does not have these properties. There, the whole mesh hierarchy needs to be recomputed with the extra added and/or exchanged level(s). The previously computed samples cannot be reused. In addition to these properties, the Local Nested Approach also has a smaller runtime for a given tolerance on the RMSE than the Global Nested Approach. Based on the results in this work, we conclude that for selecting  the random field evaluation points, an approach where the points are nested across the mesh hierarchy provides superior results compared to a non-nested approach. Of the nested approaches we consider in this paper, the Local Nested Approach has the smallest computational runtime, outperforming the Non-Nested Approach by a factor 5.

\section*{Acknowledgments}
The authors gratefully acknowledge the support from the Research Council of KU Leuven through  project C16/17/008 ``Efficient methods for large-scale PDE-constrained optimization in the presence of uncertainty and complex technological constraints". The computational resources and services used in this work were provided by the VSC (Flemish Supercomputer Center), funded by the Research Foundation - Flanders (FWO) and the Flemish Government – department EWI. 
\bibliography{Bib_ref.bib}

\end{document}